\documentclass[12pt]{article}
\usepackage{amsfonts,amsmath}

\def\be{\begin{equation}}
\def\eqn#1{\be\label{#1}}

\def\ee{\end{equation}}

\def\bea{\begin{eqnarray}}
\def\eqnn#1{\bea\label{#1}}
\def\eea{\end{eqnarray}}

\def\D{\Delta}

\newcommand{\CC}{{\mathbb C}}

\normalbaselines
\begin{document}
\pagestyle{empty}

\begin{center}
  \textsf{\LARGE
  On a "New" Deformation of GL(2)
  }

\vspace{10mm}

 {\large A.~Chakrabarti$^{a,}$\footnote{chakra@cpht.polytechnique.fr}, %\\[2mm]
 V.K.~Dobrev$^{b,c,}$\footnote{dobrev@inrne.bas.bg}
 ~and~ S.G.~Mihov$^{b,}$\footnote{smikhov@inrne.bas.bg}}

 \vspace{5mm}

 \emph{$^a$ Centre de Physique Th{\'e}orique, CNRS UMR 7644}
 \\
 \emph{Ecole Polytechnique, 91128 Palaiseau Cedex, France.}
 \\
 \vspace{3mm}
 \emph{$^b$ Institute of Nuclear Research and Nuclear Energy}
 \\
 \emph{Bulgarian Academy of Sciences}
 \\
 \emph{72 Tsarigradsko Chaussee, 1784 Sofia, Bulgaria}
 \\
 \vspace{3mm}
\emph{$^c$ Abdus Salam International Center for Theoretical Physics}
\\
\emph{Strada Costiera 11, P.O. Box 586}
\\
\emph{34100 Trieste, Italy}

\end{center}

\vspace{.8 cm}

\begin{abstract}
We refute a recent claim in the literature of a "new"
quantum deformation of GL(2).
\end{abstract}

\rightline{INRNE-TH-06-22}
%\rightline{IC/06/NN}
%\rightline{math.QA/0006206}
\rightline{December 2006}

\pagestyle{plain}
\setcounter{page}{1}

\section*{}
\label{sect:intro} \setcounter{equation}{0}

\hskip 10mm
Until the year 2000  it was not clear
how many distinct quantum group deformations are admissible for the
group $\ GL(2)\ $ and the supergroup $\ GL(1|1)\ $. For the group $\
GL(2)\ $ there were the well-known standard $GL_{pq}(2)$ \cite{DMMZ}
and nonstandard (Jordanian) $GL_{gh}(2)$ \cite{Ag} two-parameter
deformations. (The dual quantum algebras of $GL_{pq}$ and $GL_{gh}$
were found in \cite{DobrevJMP33} and \cite{ADMjpa30}, respectively.)
For the supergroup $\ GL(1|1)\ $ there were the
standard $GL_{pq}(1|1)$ \cite{HiRi,DaWa,BuTo} and the hybrid
(standard-nonstandard) $GL_{qh}(1|1)$ \cite{FHR} two-parameter
deformations.

Then, in the year 2000 in \cite{AACDM} it
was shown that the list of these four deformations is exhaustive
(refuting a long standing claim of \cite{Kup} (supported also in
\cite{BHP,LMO}) for the existence of a hybrid (standard-nonstandard)
two-parameter deformation of $GL(2)$). In particular, it was shown
that the above four deformations match the distinct triangular $4\times
4$ $R$-matrices from the classification of \cite{Hietarinta} which
are deformations of the trivial $R$-matrix (corresponding to
undeformed $GL(2)$).\footnote{Superficially, there are seven
 triangular $4\times4$ $R$-matrices in  \cite{Hietarinta},
 however, three of them are special cases of the essential four,
 cf. \cite{AACDM}}

At the end of the Introduction of \cite{AACDM} one can read the following:\\
"Instead of briefly stating the equivalence of the hybrid (standard-nonstandard)
of \cite{Kup} with the standard $GL_q(2)$, we have chosen to present our elementary
analysis explicitly and in some detail. We consider this worthwhile
for dissipating some confusions. Several authors have presented
attractive looking hybrid deformations without noticing disguised
equivalences. We ourselves devoted time and effort to their study
before reducing them to usual deformations. We hope that our analysis
will create a more acute awareness of traps in this domain."

In spite of this there still appear statements about "new" deformations of $GL(2)$.
In particular, in the Conclusions of the paper \cite{IlLy}  we read:\\ "Thus, we have a new
quantization of $GL(2)$ that is neither a twist deformation nor a quasitriangular one."

Unfortunately, the authors of \cite{IlLy} have not
noticed that their "new" quantization of $GL(2)$
is actually a partial case of the two-parameter
nonstandard (Jordanian) $GL_{gh}(2)$ deformation \cite{Ag}.

It is easy to demonstrate this explicitly. First we repeat the relations for the four generators
~$a,b,c,d,$~ of deformed  $GL(2)$ from the paper \cite{IlLy}:\\
The co-product is standard:
\eqnn{copr}
\D (a) ~=~ a \otimes a + b \otimes c \cr
\D (b) ~=~ a \otimes b + b \otimes d \cr
\D (c) ~=~ c \otimes a + d \otimes c \cr
\D (d) ~=~ c \otimes b + d \otimes d
\eea
while the algebra relations given in (10) of \cite{IlLy} are:
\eqnn{glh}
&& [a,b] ~=~ b^2 ~, \qquad [a,c] = 0 ~,\cr
&& [b,c] ~=~ -db ~, \qquad [b,d] ~=~ 0 ~, \cr
&& [a,d] ~=~ db ~, \qquad [c,d] ~=~ d^2 - ad +cb  ~.
\eea

On the other hand  the two-parameter
nonstandard (Jordanian) $GL_{gh}(2)$ deformation \cite{Ag} is given as follows.
The co-product is the standard one given above in \eqref{copr},
while the algebra relations are ($g,h \in \CC$)~:
\eqnn{glgh}
&&[d,c] ~=~ hc^2 ~, \qquad [d,b] ~=~ g(ad-bc+hac - d^2) ~, \cr
&&[b,c] ~=~ gdc+hac - ghc^2 ~, \qquad [a,c] ~=~ gc^2 ~, \cr
&& [a,d] ~=~ gdc-hac ~, \qquad [a,b] ~=~ h(da-bc+ gdc -a^2 ) ~.
\eea

It is easy to notice that \eqref{glh} is a special case of \eqref{glgh}
obtained for ~$g=0$, ~$h=1$.\\ To show this, as a first step, we set
the latter values in \eqref{glgh} to obtain:
\eqnn{glgha}
&&[d,c] ~=~ c^2 ~, \qquad [d,b] ~=~ 0  ~, \cr
&&[b,c] ~=~ ac  ~, \qquad [a,c] ~=~ 0 ~, \cr
&& [a,d] ~=~ -ac ~, \qquad [a,b] ~=~ da - bc -a^2  ~.
\eea

Now we  note that under the exchange:
\eqn{exx}  a \longleftrightarrow d
~, \qquad b \longleftrightarrow c
\ee
the co-product \eqref{copr} remains unchanged, while \eqref{glgh}
becomes:
\eqnn{glghb}
&&[a,b] ~=~ b^2 ~, \qquad [a,c] ~=~ 0  ~, \cr
&&[c,b] ~=~ db  ~, \qquad [d,b] ~=~ 0 ~, \cr
&& [d,a] ~=~ -db ~, \qquad [d,c] ~=~ ad-cb  -d^2 ~.
\eea
Clearly, \eqref{glghb} coincides with \eqref{glh}.

Thus, as anticipated there is no new deformation of $GL(2)$ in \cite{IlLy}.

\paragraph{Acknowledgments:}

%\section*{Acknowledgments}
VKD and SGM acknowledge partial support by the Bulgarian National
Council for Scientific Research, grant F-1205/02, and the European
RTN 'Forces-Universe', contract MRTN-CT-2004-005104. VKD
acknowledges partial support by the Alexander von Humboldt
Foundation in the framework of the Clausthal-Leipzig-Sofia
Cooperation.

\end{document}